\documentclass[11 pt]{article}
 \usepackage{amssymb}
\usepackage{verbatim}
\usepackage{color}

\usepackage{latexsym}
\usepackage{eucal}
\newtheorem{lemma}{Lemma}[section]
\newtheorem{proposition}{Proposition}[section]
\newtheorem{remark}{Remark}[section]

\newtheorem{sub}{\name}[section]
\newcommand{\bs}{
\begin{sub}}
\newcommand{\es}{
\end{sub}}
\newcommand{\bsl}[1]{
\begin{sub}\label{#1}}
\newcommand{\bth}[1]{\def\name{Theorem}
\begin{sub}\label{t:#1}}
\newcommand{\blemma}[1]{\def\name{Lemma}
\begin{sub}\label{l:#1}}
\newcommand{\bcor}[1]{\def\name{Corollary}
\begin{sub}\label{c:#1}}
\newcommand{\bdef}[1]{\def\name{Definition}
\begin{sub}\label{d:#1}}
\newcommand{\bprop}[1]{\def\name{Proposition}
\begin{sub}\label{p:#1}}
\newcommand{\brem}[1]{\def\name{Remark}
\begin{sub}\label{r:#1}}

\newcommand{\rth}[1]{Theorem~\ref{t:#1}}
\newcommand{\rlemma}[1]{Lemma~\ref{l:#1}}


\newcommand{\BA}{
\begin{array}}
\newcommand{\EA}{
\end{array}}
\newcommand{\bel}[1]{\begin{equation}\label{#1}}
\newcommand{\ee}{\end{equation}}%% This macro does not work with amstex.
\newcommand{\BAN}{\renewcommand{\arraystretch}{1.2}
\setlength{\arraycolsep}{2pt}
\begin{array}}
\newcommand{\BAV}[2]{\renewcommand{\arraystretch}{#1}
\setlength{\arraycolsep}{#2}
\begin{array}}
\newcommand{\BSA}{
\begin{subarray}}
\newcommand{\ESA}{
\end{subarray}}
\newcommand{\BAL}{
\begin{aligned}}
\newcommand{\EAL}{
\end{aligned}}
\newcommand{\BALG}{
\begin{alignat}}
\newcommand{\EALG}{
\end{alignat}}%% the abbrev. does not work with latex2e
\newcommand{\BALGN}{
\begin{alignat*}}
\newcommand{\EALGN}{
\end{alignat*}}%% the abbrev. does not work with latex2e
\newcommand{\note}[1]{\textit{#1.}\hspace{2mm}}
\newcommand{\Proof}{\note{Proof}}
\newcommand{\qeda}{\hspace{10mm}\hfill $\square$}

\newcommand{\Remark}{\note{Remark}}
\newcommand{\abs}[1]{\left |#1\right |}%% adjustable vertical 
\newcommand{\norm}[1]{\left \|#1\right \|}%% adjustable norm
\def\angb<#1>{\langle #1 \rangle}%% angle bracket
\newcommand{\myfrac}[2]{{\displaystyle \frac{#1}{#2} }}
\newcommand{\myint}[2]{{\displaystyle \int_{#1}^{#2}}}

\newcommand{\prt}{
\partial}

\newcommand{\ti}{\times}

\def\ga{\alpha}     \def\gb{\beta}       
             \def\ge{\epsilon}
\def\gth{\theta}                         
\def\gf{\phi}           
            \def\gl{\lambda}
\def\gm{\mu}        \def\gn{\nu}         
            
\def \gs{\sigma}       
      \def\gw{\omega}
                
     \def\Gd{\Delta}      \def\Gf{\Phi}

\def\Gl{\Lambda}          
\def\Gw{\Omega}              

   \def\CB{{\mathcal B}}   \def\CC{{\mathcal C}}

   \def\BBN {\mathbb N}    
   \def\BBR {\mathbb R}

\catcode`\@=11

\def\la{\lambda}
\def\de{\delta}

\def\eqalign#1{\null\,\vcenter{\openup1\jot \m@th
  \ialign{\strut\hfil$\displaystyle{##}$&$\displaystyle{{}##}$\hfil
     &&\strut$\displaystyle{##}$\hfil&$\displaystyle{{}##}$
     \hfil\crcr#1\crcr}}\,}
          \def\be{\begin{equation}}
     \def\ee{\end{equation}}
     \def\m{\noalign{\medskip}}
\newcommand{\rife}[1]{(\ref{#1})}

\begin{document}
\title {\bf Separable solutions of quasilinear \\Lane-Emden equations}
\author{{\bf\large Alessio Porretta}\\
 {\it Dipartimento di Matematica}\\ {\it Universit\`a di Roma Tor Vergata, Roma}\\[2mm]
{\bf\large Laurent V\'eron}\\
{\it Laboratoire de Math\'ematiques et Physique Th\'eorique}\\ {\it Universit\'e Fran\c cois Rabelais, Tours}}

\date{}
\maketitle
\noindent{\small {\bf Abstract} For $0<p-1<q$ and either $\ge=1$ or $\ge=-1$, we prove the existence of solutions of $-\Gd_pu=\ge u^q$ in a cone $C_S$, with vertex $0$ and opening $S$, vanishing on $\prt C_S$, under the form $u(x)=|x|^{-\gb}\gw(\frac{x}{|x|})$. The problem reduces to a quasilinear elliptic equation on $S$ and existence is based upon degree theory and homotopy methods. We also obtain a non-existence result in some critical case by an integral type identity.
}\smallskip

\vspace{1mm}
\hspace{.05in}

\section {Introduction}
\setcounter{equation}{0}

It is well established that the description of the boundary behavior of positive singular solutions of Lane-Emden equations
\begin{equation}\label{LE0}
-\Gd u=\ge u^q
\end{equation}
with $q>1$ in a domain $\Gw\subset \BBR^N$ is greatly helped by using specific separable solutions of the same equation. This was performed in 1991 by Gmira-V\'eron \cite{GmVe} in the case $\ge=-1$ and more recently by Bidaut-V\'eron-Ponce-V\'eron \cite{BVPV} in the case $\ge=1$. If the domain is assumed to be a cone $C_{S}=\{x\in\BBR^{N}\setminus\{0\}:x/|x|\in S\}$ with vertex $0$ and opening $S\subsetneq S^{N-1}$ (the unit sphere in $\BBR^N$), separable solutions of $(\ref{LE0} )$ vanishing on $\prt C_{_{S}}\setminus \{0\}$ were of the form
\begin{equation}\label{LE0'}u(x)=|x|^{-\frac{2}{q-1}}\gw(x/|x|),
\end{equation}  
with $\gw$ satisfying
\begin{equation}\label{LE0''}
-\Gd'\gw-\ell_{q,N}\gw-\ge\gw^q=0\qquad\mbox{in }\; S,
\end{equation}  
vanishing on $\prt S$ and where $\ell_{q,N}=\left(\left(\frac{2}{q-1}\right)\left(\frac{2q}{q-1}-N\right)\right)$ and $\Gd'$ is the Laplace-Beltrami operator on $S^{N-1}$.  To this equation is associated the functional
\begin{equation}\label{LE0'''}
J(\gf):=\myint{S}{}\left(\frac{1}{2}|\nabla'\gf|^2-\frac{\ell_{q,N}}{2}\gf^2-\frac{\ge}{q+1}|\gf|^{q+1}\right)dv_g,
\end{equation}  
where $\nabla'$ is the covariant derivative on $S^{N-1}$.
In the case $\ge=1$, non-existence of a non-trivial positive solution of $(\ref{LE0''})$ when $\ell_{q,N}\geq  \gl_{_S}$ (the first eigenvalue of $-\Gd'$ in $W^{1,2}_0(S)$) follows by multiplying the equation by the first eigenfunction and integrating over $S$; existence holds when $\ell_{q,N}< \gl_{_S}$ and $q<\frac{N+1}{N-3}$ by classical variational methods, and again non-existence holds when $q\geq \frac{N+1}{N-3}$ and $S\subset S^{N-1}_+$ is starshaped by using an integral identity \cite[Th 2.1,Cor 2.1]{BVPV}.
When $\ge=-1$, non-existence of a non-trivial solution of $(\ref{LE0''})$ when $\ell_{q,N}\leq \gl_{_S}$ is obtained by multiplying the equation by $\gw$ and integrating over $S$, while existence when $\ell_{q,N}> \gl_{_S}$ follows by minimizing $J$ over $W^{1,2}_0(S)\cap L^{q+1}(S)$.\smallskip

 In this paper we investigate similar questions for  the quasilinear  Lane-Emden equations
\begin{equation}\label{LE1}
-div\left(|\nabla u|^{p-2}\nabla u\right)=\ge u^q\qquad\mbox {in }\;C_S,
\end{equation}
where $S$ is a smooth subset of $S^{N-1}$, $q>p-1>0$ and $\ge=\pm 1$
and we look for positive solutions $u$,
vanishing on $\prt C_{_{S}}\setminus \{0\}$,  under the separable form
\begin{equation}\label{LE2}u(x)=|x|^{-\gb}\gw(x/|x|).
\end{equation}  
It is straightforward to check that $u$ is a  solution of \rife{LE1} provided 
\begin{equation}\label{LE-3}\beta=\gb_q:=\myfrac{p}{q+1-p}
\end{equation}
and $\gw$ is a positive solution of
\begin{equation}\label{E1}
-div\left(\left(\gb_q^2\gw^2+|\nabla'\gw|^2\right)^{(p-2)/2}\nabla'\gw\right)
-\gb_q\gl(\gb_q)\left(\gb_q^2\gw^2+|\nabla'\gw|^2\right)^{(p-2)/2}\gw
=\ge \gw^q
\end{equation}
in $S$ vanishing on $\prt S$, where $div(\cdot)$ is the divergence operator defined according to  the intrinsic metric $g$   and where we have set
\be\label{lab}
\gl(\gb)=\gb(p-1)+p-N.
\ee
If $\ge=0$, it is now well-known that  positive p-harmonic functions in $C_S$ vanishing on $\prt C_S$ exist under the form $(\ref{LE2})$, and either they are regular at $0$ and $\gb=-\tilde\gb_{_{S}}<0$, or 
they are singular and $\gb=\gb_{_{S}}>0$, where the values of $\tilde\gb_{_{S}}$, $\beta_S$ are unique. In this case $\gw=\tilde \gw_{_S}$ or $\gw_{_S}$  is a solution of 
\begin{equation}\label{E4}
-div\left(\left( \gb^2\gw^2+|\nabla'\gw|^2\right)^{(p-2)/2}\nabla'\gw\right)
- \gb\gl( \gb)\left( \gb^2\gw^2+|\nabla'\gw|^2\right)^{(p-2)/2}\gw
=0
\end{equation}
in $S$, where $\gb=\tilde\gb_{_{S}}$ or $\gb_{_{S}}$. The existence of $(\tilde\gb_{_{S}},\tilde \gw_{_S})$ is due to Tolksdorf in a pioneering  work \cite{To}. Tolksdorf's method has been adapted by V\'eron \cite{Ve1} in order to prove the existence of $(\gb_{_{S}},\gw_{_S})$. Later on Porretta and V\'eron \cite {PoVe} obtained a more general proof of the existence of such couples. Notice that $\gb_S$ (as well as $\tilde\gb_{_{S}}$) is uniquely determined while $\gw$ is unique up to homothety. 
In both cases the proofs rely on strong maximum principle. 

When $p\neq 2$, existence of a nontrivial solution in the case $\ge=1$ is obtained in \cite{BVJV} when $N=2$ and $\gb_{q}< \gb_{_S}$ by a dynamical system approach; while if $\ge=-1$ and $\gb_{q}> \gb_{_S}$, such an existence is proved  in  \cite{Ve1} by a suitable adaptation of Tolksdorf's construction. Notice that no functional can be associated to $(\ref{E1})$, excepted in the case  $q=q^*=\frac{ Np}{N-p} -1$. In such a  case  $(\ref{E1})$ is the Euler-Lagrange equation for the functional
\begin{equation}\label{E4'}
J_q(\gf):=\myint{S}{}\left(\frac{1}{p}\left(\gb^2_{q^*}\gf^2+|\nabla'\gf|^2\right)^{\frac{p}{2}}-\frac{\ge}{q^*+1}|\gf|^{q^*+1}\right)dv_g,
\end{equation}
and existence of a non-trivial solution of $(\ref{E1})$ with $\ge=1$ is derived from the mountain pass theorem. In all the other cases variational techniques cannot be used and have to be replaced by topological methods based upon Leray-Schauder degree. 
Define $q_c$ by
$$
q_c=q_{c,p}=\left\{\BA {ll}\frac{ (N-1)p}{N-1-p} -1& \hbox{if }p<N-1\\
\infty & \hbox{if }p\geq N-1,
\EA\right.
$$ 
then we prove the following results: \medskip

\noindent {\bf I} {\it Let $\ge=1$. Assume $p>1$, $q<q_c$ and $\gb_{q}< \gb_{_S}$, then $(\ref{E1})$  admits a positive solution in $S$ vanishing on $\prt S$}.\medskip

\noindent {\bf II} {\it Let $\ge=-1$.  Assume $p>1$ and $\gb_{q}> \gb_{_S}$, then $(\ref{E1})$ admits a unique positive solution in $S$ vanishing on $\prt S$}.\medskip

The result {\bf I} is based upon sharp Liouville theorems for solutions of $(\ref{LE1})$ in $\BBR^N$ or $\BBR^N_+$ respectively due to Serrin-Zou \cite{SeZo} and Zou \cite{Zou}. In the case of {\bf II}, the existence part is already known, but we give here a simpler form than the one in  \cite{Ve1}, using a topological deformation acting on the exponent $p$. In the case $\ge=1$, the result is optimal in the case $q=q_c$; indeed, using an integral identity, we also prove 
\medskip

\noindent {\bf III} {\it Let  $\ge=1$, $S\subsetneq S^{N-1}_+$ be a starshaped domain and  $1<p<N-1$.  If 
%$N-1>p\geq N-2$ and $q=  q_c$
%, or if $p<N-2$ and $q>p^2+p-1$, 
$q=q_c$, 
then $(\ref{E1})$ admits no positive solution in $S$ vanishing on $\prt S$}.\medskip

Notice that when $p=2$ an integral identity was used in \cite{BVPV} to prove non existence for all $q\geq q_{c,2}$. The form which is derived in the case $p\neq 2$ is much more complicated and we prove non-existence only  in the case $q= q_{c,p}$. 

Finally, the constraint $\beta_q<\beta_S$ in {\bf I} (respectivey, $\beta_q>\beta_S$ in {\bf II}) is sharp.
When $\ge= 1$, the non-existence of positive solutions of (\ref{E1}) when $\gb_{q}\geq \gb_{_S}$ has been proved  in \cite {BVJV}. The method is based upon strong maximum principle. When $\ge= -1$ a somewhat similar method is used in \cite{Ve3} and yields to non-existence results when $\beta_q\leq \beta_S$. Notice that the obtention of such results when $p=2$ is straightforward.
%%%%%%%%%%%%%%%%%%%%%%%%%%%%%%%%%%%%%%%%%%%%%%%%%%%%%%%%%%%%%%%%%%%%%%%%%%%%%%%%%%%%%%%%%%%%%%%%%%%%%%%%%%%%%%%%%%%%%%%%%%%%%%%%%%%%%%%%%%NEW-SECTION%%%%%%%%%%%%%%%%%%%%%%%%%%%%%%%%%%%%%%%%%%%%%%%%%%%%%%%%%%%%%%%%%%%%%%%%%%%%%%%%%%%%%%%%%%%%%%%%%%%%%%%%%%%%%%%%%%%%%%%%%%%%%%%%%%%%%%%%%%%%%%%%%%%%%%%%%%%%%%%%%%%%%%%%%%%%%%%%%%%%%%%%%%%%%%%%%%%%%%%%%%%%%%%%%%%%%%%%%%%%
\section {Nonexistence for the reaction problem}
\setcounter{equation}{0} 

Let $S$ be a bounded $C^2$ sub-domain of $S^{N-1}$. We consider the positive solutions in $S$ of
\be\label{pbbe}
-div\left(\left(\gb^2\gw^2+|\nabla'\gw|^2\right)^{(p-2)/2}\nabla'\gw\right)
-\gb\gl(\gb)\left(\gb^2\gw^2+|\nabla'\gw|^2\right)^{(p-2)/2}\gw
= \gw^q
\ee
vanishing on $\prt S$. Recall that $\la(\beta)$ is given by \rife{lab} and that, in connection with problem \rife{LE1}, we have interest in the special case where $\beta=\beta_q$ is given by \rife{LE-3}. The following Pohozaev-type identity, which is valid for any $\beta$,  is the key for non-existence. We denote by $S^{N-1}_+$ the half sphere. 

\begin{proposition}\label{lPPr} Let $S\subsetneq S^{N-1}$ be a $C^2$ domain and $\gf$ the first eigenfunction of $-\Gd'$ in $W^{1,2}_0(S^{N-1}_+)$. If $\omega\in W^{1,p}_0(S)\cap C(\overline S)$ is a positive solution in $S$ of \rife{pbbe}, and if  we set $\Gw=(\gb^2\gw^2+|\nabla'\gw|^2)^{1/2}$, then the  following identity holds
\begin{equation}\label{IP1}\left(1-\myfrac{1}{p}\right)\myint{\prt S}{}|\gw_{\gn}|^p\gf_{\gn}dS
=
 A\myint{S}{}\gw^{q+1}\gf\,d\gs+B\myint{S}{}\Gw^{p-2}|\nabla'\gw|^2\gf\,d\gs+C\myint{S}{}\Gw^{p-2}\gw^2\gf \,d\gs,
\end{equation}
with
\begin{equation}\label{IP2}A=A(\gb):=-\myfrac{N-1}{q+1}
-\gb(p\gb+p-N)
%=-\gb_q\left(\myfrac{N-1}{p(\gb_q+1)}+p\gb_q+p-N\right),
\end{equation}
\begin{equation}\label{IP3}B=B(\gb):=\myfrac{N-1-p}{p}+\gb(p\gb+p-N),
\end{equation}
\begin{equation}\label{IP4}C=C(\gb):=\gb^2\left(\myfrac{N-1}{p}-(p\gb+p-N)\la(\beta)\right)
%\gb_q^2\left(\myfrac{N-1}{p}-(p\gb_q+p-N)((p-1)\gb_q+p-N)\right)
.\end{equation}
\end{proposition}

In order to prove   Proposition \ref{lPPr}, we start with the following  lemma.

\begin{lemma}\label{int-part} Let $S\subset S^{N-1}$ be a $C^2$ domain and $\gf\in C^2(\overline S)$. If $\omega\in W^{1,p}_0(S)\cap C(\overline S)$ is a positive solution of  \rife{pbbe}  in $S$, we have:
%the equation 
%\be\label{pbbe}
%-div\left(\left(\gb^2\gw^2+|\nabla'\gw|^2\right)^{(p-2)/2}\nabla'\gw\right)
%-\gb\gl(\gb)\left(\gb^2\gw^2+|\nabla'\gw|^2\right)^{(p-2)/2}\gw
%=\ge \gw^q
%\ee

\begin{equation}\label{allfi}
\BA {l}\left(1-\myfrac{1}{p}\right)\myint{\prt S}{}|\gw_{\gn}|^p\gf_{\gn}dS
=
\myint{S}{}\left(\myfrac{\Gd'\gf}{q+1}
-\gb(p\gb+p-N)\gf\right) \gw^{q+1}\,d\gs-\myfrac{1}{p}\myint{S}{}\Gw^p\Gd'\gf \,d\gs\\
\m
\phantom{|\gw_{\gn}|^pdS=}
+\myint{S}{}\Gw^{p-2}D^2\gf(\nabla'\gw,\nabla'\gw) d\gs
+\gb(p\gb+p-N)\myint{S}{}\Gw^{p-2}|\nabla'\gw|^2\gf\,
d\gs\\
\m
\phantom{|\gw_{\gn}|^pdS=}
-\gb^2(p\gb+p-N)\lambda(\beta)\myint{S}{}\Gw^{p-2}\gw^2\gf\,d\sigma.
\EA
\end{equation}
\end{lemma}
%%PROOOF%%%%%%%%%%%%%%%%%%%%%
\noindent\Proof By the regularity theory of $p$-Laplace type equations (see  e.g. \cite{Dib}, \cite{To1} and the Appendix in \cite{PoVe}) it turns out that $\omega\in C^{1,\gamma}(\overline S)$ for some $\gamma\in (0,1)$, and since $\left(\gb^2\gw^2+|\nabla'\gw|^2\right)>0$ in the interior, by elliptic regularity we have $\omega\in C^2(S)$. Let   $\phi\in C^2(S)$ be a given function and $\zeta\in C^1_c(S)$; since $\zeta$ is compactly supported we can multiply \rife{pbbe} by the test function $\langle\nabla'\gw,\nabla'\gf\rangle\zeta$.   Integrating by parts we get (using the notation $\Omega:= (\gb^2\gw^2+|\nabla'\gw|^2)^{1/2}$)
$$
%\begin{equation}\label{IP5}
\BA {l}
\myint{S}{}\Gw^{p-2}\left(\myfrac{1}{2}\langle\nabla'|\nabla'\gw|^2,\nabla'\gf\rangle+D^2\gf(\nabla'\gw,\nabla'\gw)\right)\,\zeta\,d\gs
+ \myint{S}{}\Gw^{p-2}\langle\nabla'\gw,\nabla \zeta\rangle\, \langle\nabla'\gw,\nabla'\gf\rangle\,d\sigma
\\\phantom{----}
=\gb\gl(\gb)\myint{S}{}\Gw^{p-2}\gw\langle\nabla'\gw,\nabla'\gf\rangle\,\zeta\,d\gs+\myfrac{1}{q+1}\myint{S}{}\langle\nabla'\gw^{q+1},\nabla'\gf\rangle\,\zeta\,d\gs.
\EA
$$
Since
$$
\Gw^{p-2}\,\myfrac{1}{2}\langle\nabla'|\nabla'\gw|^2,\nabla'\gf\rangle= \frac1p \langle\nabla'\Omega^p,\nabla \phi\rangle-\beta^2\Omega^{p-2}\omega\langle\nabla'\gw,\nabla'\gf\rangle
$$
we obtain, due to \rife{lab}, 
$$
\BA {l}
\myfrac1p\myint{S}{} \langle\nabla'\Omega^p,\nabla' \phi\rangle\,\zeta\,d\sigma+\myint{S}{} \Gw^{p-2}\,D^2\gf(\nabla'\gw,\nabla'\gw)\,\zeta\,d\gs
+ \myint{S}{}\Gw^{p-2}\langle\nabla'\gw,\nabla' \zeta\rangle\, \langle\nabla'\gw,\nabla'\gf\rangle\,d\sigma
\\\phantom{----}
=\gb(p\beta+p-N)\myint{S}{}\Gw^{p-2}\gw\langle\nabla'\gw,\nabla'\gf\rangle\,\zeta\,d\gs+\myfrac{1}{q+1}\myint{S}{}\langle\nabla'\gw^{q+1},\nabla'\gf\rangle\,\zeta\,d\gs.
\EA
$$
Integrating by parts the first and last term we get 
\be\label{predelta}
\BA {l}
-\myfrac1p\myint{S}{}  \Omega^p  \langle\nabla' \phi,\nabla' \zeta\rangle \,d\sigma
+\myfrac{1}{q+1}\myint{S}{} \gw^{q+1}\langle\nabla'\gf,\nabla'\zeta\rangle \,d\gs
+\myint{S}{} \left(\myfrac{\gw^{q+1}}{q+1}-\myfrac{\Gw^p}{p}\right)\Gd'\gf\,\zeta \,d\gs 
\\
\phantom{----}  +\myint{S}{} \Gw^{p-2}\, D^2\gf(\nabla'\gw,\nabla'\gw)\,\zeta\,d\gs
+ \myint{S}{}\Gw^{p-2}\langle\nabla'\gw,\nabla' \zeta\rangle\, \langle\nabla'\gw,\nabla'\gf\rangle\,d\sigma \\\phantom{----} =\gb(p\beta+p-N)\myint{S}{}\Gw^{p-2}\gw\langle\nabla'\gw,\nabla'\gf\rangle\,\zeta\,d\gs.
\EA
\ee
Now we choose $\zeta= \zeta_\de$, where $\zeta_\de$ is a sequence of $C^1$ compactly supported functions such that $\zeta_\de(\sigma) \to 1$ for every $\sigma\in S$ and $|\nabla'\zeta_\de|$ is bounded in $L^1(S)$. It is easy to see by integration by parts that we have for every continuous vector field $F\in C(\overline S)$
$$
\myint{S}{} \langle F,\nabla'\zeta_\de\rangle\,d\sigma \to  - \myint{\prt S}{} \langle F,\nu(\sigma)\rangle\,d\sigma
$$
where $\nu$ is the outward unit normal on $\prt S$. We take $\zeta=\zeta_\de$ in \rife{predelta} and we let $\de\to 0$. Using that $\omega\in C^1(\overline S)$ and that, by Hopf lemma, $\omega_\nu:=\langle\nabla'\omega,\nu(\sigma)\rangle<0$ we can actually pass to the limit in the integrals containing $\nabla' \zeta_\de$. Recalling that $\omega=0$ and  $\nabla'\omega=- |\omega_\nu|  \nu$ on $\prt S$ we obtain
\begin{equation}\label{ax}
\BA {l}\left(1-\myfrac{1}{p}\right)\myint{\prt S}{}|\gw_{\gn}|^p\gf_{\gn}dS
=
\myint{S}{}\left(\myfrac{\gw^{q+1}}{q+1}-\myfrac{\Gw^p}{p}\right)\Gd'\gf \,d\gs+\myint{S}{}\Gw^{p-2}D^2\gf(\nabla'\gw,\nabla'\gw) d\gs
\\[4mm]\phantom{\left(1-\myfrac{1}{p}\right)\myint{\prt S}{}|\gw_{\gn}|^p\gf_{\gn}dS=}
-\gb(p\gb+p-N)\myint{S}{}\Gw^{p-2}\gw\langle\nabla'\gw,\nabla'\gf \rangle
d\gs.
\EA\end{equation}
Multiplying (\ref{pbbe}) by $\gw \phi$  we derive
$$
\BA {l}\myint{S}{}\Gw^{p-2}\gw\langle\nabla'\gw,
\nabla'\gf\rangle d\gs=-\myint{S}{}\Gw^{p-2}|\nabla'\gw|^2\gf \,d\gs
+\gb\gl(\gb)\myint{S}{}\Gw^{p-2}\gw^2\gf \,d\gs+\myint{S}{}\gw^{q+1}\gf \,d\gs,
\EA
$$
so that \rife{ax} becomes, replacing its last term, 
$$
\BA {l}\left(1-\myfrac{1}{p}\right)\myint{\prt S}{}|\gw_{\gn}|^p\gf_{\gn}dS
=
\myint{S}{}\left(\myfrac{\gw^{q+1}}{q+1}-\myfrac{\Gw^p}{p}\right)\Gd'\gf \,d\gs+\myint{S}{}\Gw^{p-2}D^2\gf(\nabla'\gw,\nabla'\gw) d\gs
\\[4mm]\phantom{\left(1-\myfrac{1}{p}\right)\myint{\prt S}{}|\gw_{\gn}|^p }
- \gb(p\gb+p-N) \myint{S}{} \gw^{q+1}\,\phi\,d\gs
+\gb(p\gb+p-N)\myint{S}{}\Gw^{p-2}|\nabla'\gw|^2\gf
\\
[4mm]\phantom{\left(1-\myfrac{1}{p}\right)\myint{\prt S}{}|\gw_{\gn}|^p }
-\gb^2(p\gb+p-N)\la(\beta)\myint{S}{}\Gw^{p-2}\gw^2\gf\,d\sigma.
\EA
$$
which is \rife{allfi}.
\qeda

\vskip1em
\noindent{\it Proof of Proposition \ref{lPPr}.}
We use Lemma \ref{int-part} choosing   in \rife{allfi} $\gf$ to be the first eigenfunction of $-\Gd'$ in $W^{1,2}_{0}(S^{N-1}_{+})$. Since $\Gd'\gf=(1-N)\gf$, $D^2\gf=-\gf g_{0}$, we get
\begin{equation}\label{E7}\BA {l}\left(1-\myfrac{1}{p}\right)\myint{\prt S}{}|\gw_{\gn}|^p\gf_{\gn}dS
=
-\myint{S}{}\left(\myfrac{N-1}{q+1}
+\gb(p\gb+p-N)\right) \gw^{q+1}\gf\,d\gs\\
\m
\phantom{\left(1-\myfrac{1}{p}\right)\myint{\prt S}{}|\gw_{\gn}|^p\gf_{\gn}dS=}
+\myfrac{N-1}{p}\myint{S}{}\Gw^p\gf \,d\gs
-\myint{S}{}\Gw^{p-2}|\nabla'\gw|^2 \,\phi\, d\gs
\\
\m
\phantom{\left(1-\myfrac{1}{p}\right)\myint{\prt S}{}|\gw_{\gn}|^p\gf_{\gn}dS=}
+\gb(p\gb+p-N)\myint{S}{}\Gw^{p-2}|\nabla'\gw|^2\gf\,
d\gs\\
\m
\phantom{\left(1-\myfrac{1}{p}\right)\myint{\prt S}{}|\gw_{\gn}|^p\gf_{\gn}dS=}
-\gb^2(p\gb+p-N)\la(\beta)\myint{S}{}\Gw^{p-2}\gw^2\gf.
\EA\end{equation}
Then, using also the definition of $\Omega$,  $(\ref{IP1})$ follows, with $A$, $B$ and $C$ given by  $(\ref{IP2})$-$(\ref{IP4})$.
\qeda\medskip

We shall say that a $C^2$ domain $S\subset S^{N-1}_+$ is starshaped if there exists a spherical harmonic $\gf$ of degree $1$ such that $\gf>0$ on $S$ and for any $a\in\prt S$,
\begin{equation}\label{*Sh}\langle\nabla\gf,\gn_a\rangle\leq 0
\end{equation}
where $\gn_a$ is the unit outward normal vector to $\prt S$ at $a$ in the tangent plane $T_a$ to $S^{N-1}$. It also means  that there exists some $x_0\in S$ such that the geodesic  connecting $x_0$ and $a$ remains inside $S$.

\bth{IPTh} Assume that $1<p<N-1$, $q= q_c$ and $S\subset S^{N-1}_+$ is starshaped. Then $(\ref{pbbe})$ admits no positive solution in $S$ vanishing on $\prt S$. 
%Furthermore, the previous non-existence statement still holds \smallskip

%\noindent (i) either if $p\leq 2$ and $q> q_c$, \smallskip

%\noindent  (ii) or if $p>2$ and $q> $.
\es
%%%%%
\noindent\Proof  Recall that in \rife{E1} we have $\beta_q=\frac p{q-(p-1)}$, hence different values of $q$ are in one-to-one correspondence with different values of $\beta$.  We first notice that, if $q=q_c$ the corresponding critical $\beta$ is given by 
\be\label{betac}
\beta_c\,:= \frac p{q_c-(p-1)}=\frac{N-1-p}p\,.
\ee
We use now Proposition \ref{lPPr} with $\beta=\beta_q$ and we analyze the values of the coefficients $A$, $B$, $C$ given by $(\ref{IP2})$-$(\ref{IP4})$ as functions of $\beta$. First of all, since $q+1=  \frac{p(1+\beta)}\beta$, we have
$$
A= -\myfrac{(N-1)\beta}{p(1+\beta)}
-\gb(p\gb+p-N)= -\frac\beta{(\beta+1)} \left(\frac{N-1}p + p(\beta+1)^2-N(\beta+1)\right)
$$
and since from \rife{betac} we have $\beta_c+1= \frac{N-1}p$ we deduce
$$
A= -\frac\beta{(\beta+1)} \, p \left(\beta+1-\frac 1p\right)(\beta-\beta_c)\,.
$$
Still using \rife{betac}, we also get
$$
B= \beta_c+ \beta(p(\beta-\beta_c)-1)= (\beta-\beta_c)(\beta p-1)\,.
$$
Finally, using  \rife{lab} and \rife{betac} we have  

\begin{equation}\label{c}\BA {l}
C=\gb^2\left(\myfrac{N-1}{p}-(p\gb+p-N)((p-1)\gb+p-N)\right)
\\[2mm]\phantom{C}
=\gb^2\left(\beta_c+1-(p(\gb-\beta_c)-1)(p(\beta-\beta_c)-(\beta+1))\right)
\\[2mm]\phantom{C}
= \beta^2(\beta-\beta_c)(1-p)\left(p\beta-1-\frac{p(N-p)}{p-1}\right).
\EA
\end{equation}

Therefore   $A\geq 0$, $B\geq 0$ and $C\geq 0$ can be obtained only if $q=q_c$, i.e. $\beta=\beta_c$, in which case $A=B=C=0$. Since $\gf_{\gn}\leq 0$ because $S$ is star--shaped, we deduce from  \rife{IP1}  that  $|\gw_{\gn}|^p\gf_{\gn}=0$ on $\prt S$. Unless $\gw$ is  identically zero, we have $\gw_{\gn}<0$ by Hopf boundary lemma.  Then $\gf_{\gn}\equiv 0$, and using the equation satisfied by $\gf$ and Gauss formula, we derive
$$
\gl_{_S}\myint{S}{}\gf d\gs=0\Longrightarrow \gf\equiv 0\qquad \hbox{in $S$},
$$
which is  impossible since $\phi>0$ in $S^{N-1}_+$. This proves the first assertion.\qeda
%\qeda
%% we consider $A$, $B$ and $C$ as function of $\gb$ over the interval  $(0,\gb_{q_c})$. Then
%%$$\BA {l}\gb^{-1}(\gb+1)^{-1}A(\gb)=-p\gb^2+(N-2p)\gb-p-\myfrac{N-1}{p}+N\\[3mm]
%%\phantom{\gb^{-1}(\gb+1)^{-1}A(\gb)}
%%=-p(\gb-\gb_{q_c})(\gb-\frac{1-p}{p}).
%%\EA$$
%%Thus 
%%\begin{equation}\label{E8}
%%\gb\in (0,\gb_{q_c}]\Longrightarrow A(\gb)\geq 0.
%%\end{equation}
%%Next
%%$$B(\gb)=p(\gb-\gb_{q_c})(\gb-p^{-1}).
%%$$
%%Then
%%\begin{equation}\label{E'9}
%%\gb\in (0,\inf\{p^{-1};\gb_{q_c}]\cup [\sup\{\gb_{q_c};p^{-1}\},\infty)\Longrightarrow B(\gb)\geq 0.
%%\end{equation}
%%Finally, 
%%$$\gb^{-2}C(\gb)=-p(p-1)(\gb-\gb_{q_c})(\gb-\tilde\gb),
%%$$
%%with
%%$$\tilde \gb=-\myfrac{N-1}{p(p-1)(N-1-p)},
%%$$
%%thus
%%\begin{equation}\label{E'10}
%%\gb\in (0,\gb_{q_c}]\Longrightarrow C(\gb)\geq 0.
%%\end{equation}
%%Since $A(\gb)$ and $C(\gb)$ are nonnegative when $\gb\in (0,\gb_{q_c}]$, the condition $(\ref{E'9})$ is satisfied if $\gb\leq \inf\{p^{-1};\gb_{q_c}\}$. When $p\geq N-2$, the infimum is $\gb_{q_c}$ and we get (i); if $p< N-2$,  the infimum is $p^{-1}$ and we get (ii).\qeda\medskip
%\noindent\Remark When $p=2$, the second statement of our result reads as follows:\smallskip

%\noindent (i) either $N=3,4$ and $q\geq q_c$,\smallskip

%\noindent (ii) or $N\geq 5$ and $q\geq 5$.\smallskip
\medskip

\noindent\Remark
If $p=2$, it is proved in \cite{BVPV} that the nonexistence result of \rth{IPTh} holds for every  $q\geq q_c$, which suggests that our result  above is not optimal. The proof in \cite{BVPV} cannot be applied here since the term 
$\int_S\Gw^{p-2}\gw\langle\nabla'\gw,\nabla'\gf \rangle d\gs$ is completely integrable only if $p=2$. However, we conjecture that, even when $p\neq 2$, the conclusion of \rth{IPTh} holds under the more general  condition $q\geq q_c$. \medskip

\noindent\Remark If we assume that $p\neq 2$, the proof of \rth{IPTh} relies on the existence of  a positive function $\gf$ in $S$, satisfying $(\ref{*Sh})$ on $\prt S$ and 
\begin{equation}\label{E9}\myfrac{\Gd'\gf}{(q+1)\gf}-\gb (p\gb+p-N)\geq 0,
\end{equation}
\begin{equation}\label{E10}\myfrac{pD^2\gf(\xi,\xi)-\Gd'\gf}{p\gf}+\gb(p\gb+p-N)\geq 0\quad\forall\xi\in S^{N-1},
\end{equation}
and
\begin{equation}\label{E11}-\myfrac{\Gd'\gf}{p\gf}-(p\gb+p-N)((p-1)\gb+p-N)\geq 0.
\end{equation}

\begin{remark}\label{Re2}  {\rm For the sake of completeness, we recall the non-existence result obtained in \cite[Th 1]{BVJV}:\smallskip

\noindent {\it Let $\ge=1$ and $0<p-1<q$. If $\gb_q\geq \gb_S$, there exists no positive solution
of $(\ref{E1})$ in $S$ which vanishes on $\prt S$.}\
}\end{remark}

%%%%%%%%%%%%%%%%%%%%%%%%%%%%%%%%%%%%%%%%%%%%%%%%%%%%%%%%%%%%%%%%%%%%%%%%%%%%%%%%%%%%%%%%%%%%%%%%%%%%%%%%%%%%%%%%%%%%%%%%%%%%%%%%%%%%%%%%%%NEW-SECTION%%%%%%%%%%%%%%%%%%%%%%%%%%%%%%%%%%%%%%%%%%%%%%%%%%%%%%%%%%%%%%%%%%%%%%%%%%%%%%%%%%%%%%%%%%%%%%%%%%%%%%%%%%%%%%%%%%%%%%%%%%%%%%%%%%%%%%%%%%%%%%%%%%%%%%%%%%%%%%%%%%%%%%%%%%%%%%%%%%%%%%%%%%%%%%%%%%%%%%%%%%%%%%%%%%%%%%%%%%%%

\section {Existence for the reaction problem}
\setcounter{equation}{0} 

Concerning the problem with reaction we consider a more general statement than Theorem I, replacing the sphere by a complete  $d$-dimensional Riemannian manifold $(M,g)$ and suppose that $S$ is a relatively compact smooth open domain of $M$. We denote by $\nabla:=\nabla_g$ the gradient of a function identified with its covariant derivatives and by $div:=div_g$ the intrinsic divergence operator acting on  vector fields. The following result is proved in $\cite{PoVe}$. 

\bth{spec}For any $\gb>0$ there exists a unique $\Gl_\gb>0$ and a unique (up to an homothety)  positive function $\gw_{\gb}\in C^2( S)\cap C^1(\overline S)$ solution of
\begin{equation}\label{F1}
\left\{\BA {l}
-div\left(\left(\gb^2\gw_{\gb}^2+|\nabla \gw_{\gb}|^2\right)^{\frac{p-2}2}\nabla \gw_{\gb}\right)
=\gb\Gl_\gb\left(\gb^2\gw_{\gb}^2+|\nabla \gw_{\gb}|^2\right)^{\frac{p-2}2}\gw_{\gb}\quad \mbox{in } S\\
\m
%\phantom{div\left(\left(\gb^2\gw_{\gb}^2+|\nabla\gw_{\gb}|^2\right)^{p-2}\nabla \gw_{\gb}\right)
%+\gb\Gl_\gb\left(\gb^2\gw_{\gb}^2+|\nabla \gw_{\gb}|^2\right)^{p-2}}
\gw_{\gb}=0\quad\mbox{on } \prt S.
\EA\right.\end{equation}
The mapping $\gb\mapsto \Gl_\gb$ is continuous and decreasing, and the spectral exponent $\gb_{_S}$ is the unique $\gb>0$ such that $\Gl_{\gb_{_S}}=\gb_{_S}(p-1)+p-d-1$.
\es

\vskip0.5em

\begin{remark}\label{moniff}
Let us notice that the monotone character of $\beta\mapsto \Lambda_\beta$  implies that
$$
0<\gb<\gb_{_S}\Longleftrightarrow\Gl_\gb-\gb(p-1)>\Gl_{\gb_{_S}}-\gb_{_S}(p-1)=p-d-1
$$
Therefore, if we set $\la(\beta)=\beta(p-1)+p-d-1$, we deduce that 
\begin{equation}\label{F2}
0<\gb<\gb_{_S}\Longleftrightarrow \Gl_\gb>\gl(\gb).
\end{equation}
\end{remark}

\vskip1em

Let us now prove the existence of solutions for the reaction problem.

\bth{Exist} Assume $1<p<d$ and $p-1<q<q_c:=pd/(d-p)-1$. Then for any $0<\gb<\gb_{_S}$, there exists 
 a positive function $\gw\in C(\overline S)\cap C^2(S)$  satisfying
 \begin{equation}\label{source}
\left\{ \BA {l}
-div\left((\gb^2\gw^2+|\nabla\gw|^2)^{\frac{p}{2}-1}\nabla \gw\right)=\gb\gl(\gb)(\gb^2\gw^2+|\nabla \gw|^2)^{\frac{p}{2}-1}\gw+ \gw^q\quad\mbox{ in }\; S\\[2mm]
%\phantom{div\left((\gb^2\gw^2+|\nabla\gw|^2)^{\frac{p}{2}-1}\nabla \right)+\gb\gl(\gb)(\gb^2\gw^2+|\nabla \gw|^2)^{\frac{p}{2}-1}\gw+ \gw^q}
\gw=0\quad\mbox{ on }\;\prt S,
 \EA\right.
 \end {equation}
 where $\gl(\gb)=\gb(p-1)+p-d-1$.
\es

In order to prove \rth{Exist}, we use topological arguments as it is often needed in a non-variational setting. In particular,   following a  strategy similar as in  \cite{QS}, our proof is based upon the following  fixed point theorem  which is only one possible consequence of Leray--Schauder degree theory to compute the fixed point index of compact mappings. Such results were developed mostly by Krasnoselskii (\cite{Kr}), we refer to Proposition 2.1 and Remark 2.1 in \cite{dFLN} for the statement below. 
%%%%%%%%%%%%TH%%%%%%%%%%%%%%%%%%%%%%%%%%
\bth{figue}
Let $X$ be a Banach space and $K\subset X$ a closed cone with non empty interior. Let $F\,:K\times \BBR_+\to K$ be a compact mapping, and let $\Phi(u)=F(u,0)$ (compact mapping from $K$ into $K$).
Assume the following holds: there exist $R_1<R_2$ and $T>0$ such that
\vskip0.2em
\noindent (i) $u\neq s\Phi(u)$ for every $s\in [0,1]$ and every $u$: $\|u\|=R_1$.
\vskip0.2em
\noindent (ii) $F(u,t)\neq u$ for every $(u,t)$: $\|u\|\leq R_2$ and $t\geq T$.
\vskip0.2em
\noindent (iii) $F(u,t)\neq u$ for every $u$: $\|u\|=R_2$ and every $t\geq 0$.
\vskip0.2em
\noindent Then, the mapping $\Phi$ has a fixed point $u$ such that $R_1<\|u\|<R_2$.
\es

%%%%%TH%%%%%%%%%%%%%%%%%%%%%%%%%%%%%

We also recall the following non-existence results respectively due to Serrin and Zou \cite{SeZo}, and Zou \cite{Zou}.
\bth{SZ}
Assume $1<p<d$ and $p-1<q<q_c$. Then there exists no $C^1$ positive solution of
 \begin{equation}\label{Glob}
-\Gd_p u=u^q
 \end {equation}
 in $\BBR^d$.
\es
%%%%%%TH%%%%%%%%%%%%%%%%%%%%%%%%%%%%%%%%%%%%%%%

\bth{Z}
Assume $1<p<d$ and $p-1<q<q_c$. Then there exists no $C^1$ positive solution of
 \begin{equation}\label{Glob}
-\Gd_p u=u^q
 \end {equation}
 in $\BBR_+^d:=\{x=(x_1,..., x_d):x_d>0\}$ vanishing on $\prt\BBR_+^d:=\{x=(x_1,..., x_d):x_d=0\}$.
\es
%%%%%%%%%%%%%%%%%%%%%%%%%%%%%%%%%%%%%%%%%%%%%%%%%%%%%%%%%%%%%%%%%%%%%%%%%%%%%%%%%%%%%%%%%%%%%%%%%%%%%%%%%%%%%%%%%%%%%%%%%%%%%%%%%%%%%%%%%%%%%%%%%%%%%%%%%%%%%%%%%%%%%%%%%%%%%%%%%%%%%%OLD%%%%%%%%%%%%%%%%%%%%%%%%%%%%%%%%%%%%

\noindent{\it Proof of \rth{Exist}. }
Define the operator ${\cal A}$ in $W^{1,p}_{0}(S)$ as 
$$
{\cal A}(\gw):=-div_{g}\left((\gb^2\gw^2+|\nabla\gw|^2)^{\frac{p}{2}-1}\nabla \gw\right)+\gb^2\,\gw(\gb^2\gw^2+|\nabla \gw|^2)^{\frac{p}{2}-1}\,.
$$
Note that ${\cal A}$ is the derivative of the functional
$$
J(w)=\frac1p \int_S (\gb^2\gw^2+|\nabla\gw|^2)^{\frac{p}{2}}\,dv_g
$$
Since $J$ is strictly convex, then ${\cal A}$ is a strictly monotone  operator from $W^{1,p}_0(S)$ into $W^{-1,p'}(S)$, henceforth its inverse is well defined and continuous \cite {JLL}. 
In order to apply  \rth{figue}, we denote by  $X=C_0^1(\overline S)$, the closure of $C_0^1(S)$  in $C^1(\overline S)$. Clearly $X\subset W^{1,p}_0(S)$, with continuous imbedding, if it is endowed with its natural norm $||.||_X:=||.||_{C^1(\overline S)}$. Furthermore, since $\prt S$ is $C^2$, $C^1(\overline S)\cap W^{1,p}_0(S)=C_0^1(\overline S)$. If $K$ is the cone of nonnegative functions in $S$, it has a nonempty interior. For $t>0$, we set
$$
F(\gw,t):={\cal A}^{-1}\left( \beta\left(\la(\beta)+\beta+t\right)\gw (\gb^2\gw^2+|\nabla\gw|^2)^{\frac{p}{2}-1}+(\gw+t)^q  \right)\,.
$$
Note that 
$$
\Phi(\gw):=F(\gw,0)={\cal A}^{-1}\left( \beta\left(\la(\beta)+\beta\right)\gw (\gb^2\gw^2+|\nabla\gw|^2)^{\frac{p}{2}-1}+\gw^q  \right);
$$
henceforth any nontrivial fixed point for $\Phi$ would solve problem \rife{source}. \smallskip

We have to verify the assumptions of \rth{figue}. First of all, the compactness of $F(\gw, t)$. If we set 
$F(\gw,t)=\gf$, then it means that $\gf\in W^{1,p}_0(S)$ satisfies
\bel{eq1}\BA {l}-div_{g}\left((\gb^2\gf^2+|\nabla\gf|^2)^{\frac{p}{2}-1}\nabla \gf\right)+\gb^2\,\gf(\gb^2\gf^2+|\nabla \gf|^2)^{\frac{p}{2}-1}\\[2mm]\phantom{--------}=\left( \beta\left(\la(\beta)+\beta+t\right)\gw (\gb^2\gw^2+|\nabla\gw|^2)^{\frac{p}{2}-1}+(\gw+t)^q  \right)
\EA\ee
Thus, if we assume that $\gw$ belongs to a bounded set in $K\cap X$, the right-hand side of (\ref{eq1}) is bounded in $C(\overline S)$. Thus, by standard regularity estimates up to the boundary for $p$-Laplace type operators (see \cite[Appendix]{PoVe} and \cite{Dib}, \cite {To1}), $\gf$ remains bounded in  $C^{1,\alpha}(\overline S)$ and therefore relatively compact in $C^1(\overline S)$. It remains to show that conditions (i)--(iii) of \rth{figue} hold.
\medskip

\noindent{\it Step 1: Condition (i) holds.}  
 We proceed by contradiction in supposing that there exists a sequence $\{s_n\}\subset [0,1]$ such that for any $n\in\BBN$ the following problem
\begin{equation}\label{i-s}
\left\{ \BA {ll}
-div_{g}\left((\gb^2\gw^2+|\nabla\gw|^2)^{\frac{p}{2}-1}\nabla \gw\right)+\gb^2(\gb^2\gw^2+|\nabla \gw|^2)^{\frac{p}{2}-1}\gw\\[2mm]
\qquad\qquad\quad\quad\quad \!=s^{p-1}\gb(\gl(\gb)+\beta)(\gb^2\gw^2+|\nabla \gw|^2)^{\frac{p}{2}-1}\gw+ s_n^{p-1}\gw^q\qquad\mbox{ in }\; S\\[2mm]
\qquad\quad\quad\quad\quad 
\gw=0\quad\mbox{ on }\;\prt S,
 \EA\right.
\end{equation}
 admits a positive solution $\gw_n$, and that there holds
$$
\|\gw_n\|_{X}\to 0\quad \hbox{as $n\to \infty$}.
$$
Set $w_n=\gw_n/\|\gw_n\|$, then $w_n$ solves
$$
\left\{ \BA {l}
-div_{g}\left((\gb^2w_n^2+|\nabla w_n|^2)^{\frac{p}{2}-1}\nabla w_n\right)+\gb^2w_n(\gb^2w_n^2+|\nabla w_n|^2)^{\frac{p}{2}-1}\\
\qquad\quad =s_n^{p-1}\gb(\gl(\gb)+\beta)(\gb^2 w_n^2+|\nabla w_n|^2)^{\frac{p}{2}-1}w_n+ s_n^{p-1}w_n^q\,\|w_n\|_X^{q-(p-1)}\quad\mbox{ in }\; S
\\
\m
w_n=0\quad\mbox{ on }\;\prt S
 \EA\right.
$$
Up to subsequences, we assume that $s_n\to  s$ for some $s\in [0,1]$. Using  compactness arguments we deduce that $w_n$ will converge strongly in $C^1(\overline S)$ to some positive function $w$ such that $\|w\|_X=1$ and which solves
\begin{equation}\label{Z0}
\left\{ \BA {l}
-div_{g}\left((\gb^2w^2+|\nabla w|^2)^{\frac{p}{2}-1}\nabla w \right)
\\
\qquad\qquad = \gb\left(s^{p-1}\gl(\gb)+(s^{p-1}-1)\beta\right)(\gb^2 w^2+|\nabla w |^2)^{\frac{p}{2}-1}w\quad\mbox{ in }\; S\\
\m
%\phantom{-div_{g}\left((\gb^2\gw^2+|\nabla\gw|^2)^{\frac{p}{2}-1}\nabla\right)}
w=0\quad\mbox{ on }\;\prt S
 \EA\right.
\end{equation}
Using \rth{spec}, we derive $\Gl_\gb=s^{p-1}\gl(\gb)+(s^{p-1}-1)\beta$.
Since $\gb<\gb_{_S}$, $\la(\beta)<\Lambda_\beta$ by $(\ref{F2})$.  Therefore, as $s\leq 1$, we get
$$
s^{p-1}\gl(\gb)+(s^{p-1}-1)\beta \leq s^{p-1}\gl(\gb)<\Lambda_\beta\,,
$$
which is a contradiction. Consequently, there exists $R_1>0$ such that for any $s\in [0,1]$, there holds  $\gw\neq s\Gf (\gw)$ for any $\gw$ such that $\norm\gw_X= R_1$.\medskip

\noindent {\it Step 2: Condition (ii) holds.}  Consider the first eigenvalue $\la_{1,\gb}$ associated with the operator ${\cal A}$, i.e. 
\be\label{la1be}
\la_{1,\gb}=\min\left\{ \int_S (\gb^2\gw^2+|\nabla\gw|^2)^{\frac{p}{2}}\,dv_g: \gw\in W^{1,p}_0(S)\,,\, \int_S|\gw|^p\,dv_g=1\right\}
\ee
Note that for $t$ large enough, we have $\la(\beta)+\beta+t\geq 0$, hence, using that $q>p-1$, we can find $T>0$ such that
$$
\beta\left(\la(\beta)+\beta+t\right)\gw (\gb^2\gw^2+|\nabla\gw|^2)^{\frac{p}{2}-1}+(\gw+t)^q \geq (\la_1+\de)\gw^{p-1}\qquad\forall t\geq T\,,\forall \gw\geq 0\,.
$$
Therefore, if $t\geq T$ and $F(\gw,t)=\gw$   we deduce that $\gw\neq 0 $ and  satisfies
$$
\left\{ \BA {l}
-div_{g}\left((\gb^2\gw^2+|\nabla\gw|^2)^{\frac{p}{2}-1}\nabla \gw\right)+\gb^2\gw(\gb^2\gw^2+|\nabla \gw|^2)^{\frac{p}{2}-1}\geq (\la_{1,\gb}+\de)\gw^{p-1}\,\,\mbox{ in }\; S\\[2mm]
%\phantom{-div_{g}\left((\gb^2\gw^2+|\nabla\gw|^2)^{\frac{p}{2}-1}\nabla\right)}
\gw=0\quad\mbox{ on }\;\prt S
 \EA\right.
$$
The existence of a positive super-solution with $\la_{1,\gb}+\de$ would make it possible to construct  a  positive solution as well.
But since $\la_{1,\gb}$ is an isolated eigenvalue (see Appendix) this yields a contradiction. Therefore, for $t\geq T$ the equation $F(\gw,t)=\gw$ has no solution at all. Note that $T$ only depends on $\la_1$, $\beta$.\medskip

\noindent {\it Step 3: Condition (iii) holds.} 
Since we proved that  (ii) holds independently on the choice of $R_2$, it is enough to show that (iii) holds for every $t\leq T$. 

This is done if we have  the existence of universal a priori estimates, i.e. if we can prove  the existence of  a constant $R_2$ such that for any $t\leq T$ every positive solution of
$$
\left\{ \BA {l}
-div_{g}\left((\gb^2\gw^2+|\nabla\gw|^2)^{\frac{p}{2}-1}\nabla \gw\right)+\gb^2\gw(\gb^2\gw^2+|\nabla \gw|^2)^{\frac{p}{2}-1}=\\
\qquad\qquad  \gb(\gl(\gb)+\beta+t)(\gb^2\gw^2+|\nabla \gw|^2)^{\frac{p}{2}-1}\gw+ ( \gw+t)^q\quad\mbox{ in }\; S\\[2mm]
\phantom{-div_{g}\left((\gb^2\gw^2+|\nabla\gw|^2)^{\frac{p}{2}-1}\nabla\right)}
\gw=0\quad\mbox{ on }\;\prt S
 \EA\right.
$$
satisfies $\|\gw\|< R_2$.

The crucial step is to prove that there exist universal  a priori estimates for the $L^\infty$-norm (a  bound for the $W^{1,p}_0$-norm would follow immediately, and then  a bound in $X$ from the regularity theory). A standard procedure is to reach this result reasoning by contradiction and using a blow-up argument. Indeed, if  a universal bound does not exist, there exist a sequence of solutions $\gw_n$ and $t_n\leq T$ such that 
$$
\|\gw_n\|_\infty\to \infty\,.
$$
Let $\sigma_n$ be the (local coordinates of) maximum points of $\gw_n$; up to subsequences, we have $\sigma_n\to \sigma_0\in \overline S$. Setting $M_n=\|\gw_n\|_\infty^{-\frac{q-(p-1)}p}$, define
$$
 v_n(y)= \frac{\gw_n(\sigma_n+M_ny)}{\|\gw_n\|_\infty}= M^{\frac{p}{q-(p-1)}}_n\,\gw_n(\sigma_n+M_ny)
$$
Then $v_n$ is  a sequence of uniformly bounded solutions, which will be locally compact in the $C^1$-topology.
Rescaling the equation and passing to the limit in $n$  we find out that the limit function $v$ is positive and satisfies the equation
$$
-\Delta_p v=c_0 v^q
$$
for some constant $c_0$ (coming out from the local expression of Laplace-Beltrami operator).  Depending whether $\sigma_0\in S$ or $\sigma_0\in \partial S$, the equation would take place in either $\BBR^d$ or in the half space $\BBR^d_+$, where $d=N-1$, in which case $v$ vanishes on $\prt \BBR^d_+$. Since $p-1<q<q_c$, this contradicts either \rth {SZ}, or \rth{Z} because, by construction, we have $v(0)=1$.\qeda\medskip

\noindent\Remark In the case $p=2$, existence is proved in \cite{BVPV} using a standard variational method. It is also proved that, if $(M,g)=(S^d,g_0)$ (the standard sphere), and if $S$ is a spherical cap with center $a$, any positive solution of 
 \begin{equation}\label{source:p=2}
\left\{ \BA {l}
\Gd' \gw+\gb(\gb+1-d))\gw+ \gw^q=0\quad\mbox{ in }\; S\\[2mm]
\phantom{\Gd' \gw+\gb(\gb+1-d))+ \gw^q}
\gw=0\quad\mbox{ on }\;\prt S,
 \EA\right.
 \end {equation}
depends only on the angle $\gth$ from $a$. Furthermore, uniqueness is proved by a delicate analysis of the non-autonomous second order O.D.E. satisfied by $\gw$. In the case $p\neq 2$ and  assuming always that $S$ is a spherical cap of $(S^d,g_0)$, it is still possible to construct a radial (i.e. depending only on $\gth$) positive solution of $(\ref{source})$: it suffices to restrict the functional analysis framework to radial functions. However, there are two interesting open questions the answer to which would be important:\smallskip

\noindent (i) {\it Are all positive solutions of $(\ref{source})$ radial ?}\smallskip

\noindent (ii) {\it Is there uniqueness of positive  radial solutions of $(\ref{source})$?}

%%%%%%%%%%%%%%%%%%%%%%%%%%%%%%%%%%%%%%%%%%%%%%%%%%%%%%%%%%%%%%%%%%%%%%%%%%%%%%%%%%%%%%%%%%%%%%%%%%%%%%%%%%%%%%%%%%%%%%%%%%%%%%%%%%%%%%%%%%NEW-SECTION%%%%%%%%%%%%%%%%%%%%%%%%%%%%%%%%%%%%%%%%%%%%%%%%%%%%%%%%%%%%%%%%%%%%%%%%%%%%%%%%%%%%%%%%%%%%%%%%%%%%%%%%%%%%%%%%%%%%%%%%%%%%%%%%%%%%%%%%%%%%%%%%%%%%%%%%%%%%%%%%%%%%%%%%%%%%%%%%%%%%%%%%%%%%%%%%%%%%%%%%%%%%%%%%%%%%%%%%%%%%

%\vskip1em
\section {Existence for the absorption problem}
\setcounter{equation}{0} 
Let us now consider the absorption problem, namely \rife{E1} with $\ge=-1$.  We give an existence result which extends the previous ones obtained in \cite{Ve1}, with a  simpler proof.   

\bth{ExistII} Assume $0<p-1<q$. Then for any $\gb>\gb_{_S}$, there exists 
 a unique positive function $\gw\in C(\overline S)\cap C^2(S)$  satisfying
 \begin{equation}\label{absor}
\left\{ \BA {l}
-div_g\left((\gb^2\gw^2+|\nabla\gw|^2)^{\frac{p}{2}-1}\nabla \gw\right)=\gb\gl(\gb)(\gb^2\gw^2+|\nabla \gw|^2)^{\frac{p}{2}-1}\gw- \gw^q\quad\mbox{ in }\; S\\[2mm]
%\phantom{div_g\left((\gb^2\gw^2+|\nabla\gw|^2)^{\frac{p}{2}-1}\nabla \right)+\gb\gl(\gb)(\gb^2\gw^2+|\nabla \gw|^2)^{\frac{p}{2}-1}\gw+ \gw^q}
\gw=0\quad\mbox{ on }\;\prt S,
 \EA\right.
 \end {equation}
 where $\gl(\gb)=\gb(p-1)+p-d-1$.
\es

Before proving  \rth{ExistII}, we will need the following lemma.

\blemma{contp} For $\beta>0$ and $p>1$, let $\Lambda_\beta$ and $\beta_S$ be defined by \rth{spec}. Then both $\Lambda_\beta$ and $\beta_S$ are continuous functions of $p$, varying in $(1,\infty)$.
\es
%%%%%%%%%%%%%%%%%%%%%%%%%%%%%%%%%%%%%%%%%

\noindent\Proof   By \rth{spec},  $\Lambda_\beta$ is uniquely defined for any fixed $p>1$. To emphasize the dependence of $\Lambda_\beta$ on $p$, let us  denote it now by $\Lambda_{\beta,p}$. The continuity of $\Lambda_{\beta,p}$ with respect to $p$ can be proved in the same way  as we proved (see  Proposition 2.4 in \cite{PoVe}) the continuity of $\Lambda_{\beta,p}$ with respect  to $\beta$. Thus,  we only sketch the argument, which relies on the construction itself of $\Lambda_{\beta,p}$. Indeed, we proved in \cite{PoVe}
 that $\Lambda_{\beta,p}$ is the unique constant such that there exists a function $v\in C^2(S)$ satisfying
 \begin{equation}\label{old}\left\{ \BA {l}
 -\Gd_{g}v-(p-2)\myfrac{D^2v \nabla v.\nabla v}{1+|\nabla v|^2}+\gb(p-1)|\nabla v|^2=-\Lambda_{\beta,p}\quad\mbox{ in }S
 \\
 \m
 \lim\limits_{\sigma\to\prt S}v( \sigma)=\infty.
 \EA\right.
 \end {equation}
If we normalize $v$ by setting, for example, $v(\sigma_0)=0$ for some $\sigma_0\in S$, then $v$ is  unique. Moreover $v\in C^2(S)$ and $v$ satisfies estimates in $W^{1,\infty}_{loc}(S)$ which are uniform as $\beta\in (0,\infty)$ and $p\in (1,\infty)$ vary in  compact sets. It is also easy to check (see \cite{PoVe}) that $\Lambda_{\beta,p}$ remains bounded whenever $\beta$ varies in a compact set of $(0,\infty)$ and $p$ vary in a compact set of $(1,\infty)$.  The estimates obtained on $v$ and $\nabla v$ imply that, whenever $\beta_n$ or $p_n$ are convergent sequences, the  sequence of corresponding solutions $v_n$ of \rife{old} (such that $v_n(\sigma_0)=0$)   is relatively compact (locally uniformly in $C^1$). The equation \rife{old} turns out then to be stable (including the boundary estimates);
finally, the uniqueness property of $\Lambda_{\beta,p}$, and of the associated (normalized) solution $v$, implies the continuity of $\Lambda_{\beta,p}$ with respect to both $\beta$ and $p$.

Let now $\beta_{S,p}$ be the spectral exponent defined by the equation
\be\label{eqbetas}
\Lambda_{\beta,p}=\beta (p-1)+p-d-1
\ee
First of all note that when $p$ lies in a  compact set in $(1,\infty)$,  then necessarily $\beta_{S,p}$ is bounded. Indeed, 
since $\Lambda_{\beta,p}\leq \Lambda_{1,p}$ whenever $\beta\geq 1$, we have that  
$$
\beta_S (p-1)+p-d-1 \leq \Lambda_{1,p} \quad \hbox{if $ \beta_S\geq 1$,}
$$
so that
$$
\beta_S\leq 1 + \frac1{p-1}\left(\Lambda_{1,p}- (p-d-1)\right)\,.
$$
Therefore,  if $p$ belongs to a compact  set in $(1,\infty)$, then  $\beta_S$ remains also in  a bounded set. 
Now, if $p_n\to p_0$, setting $\beta_n=\beta_{S,p_n}$,  we have that $\beta_n$ is bounded and, up to subsequences, it is   convergent to some $\beta_0$. From \rife{eqbetas}, we deduce that $\Lambda_{\beta_n, p_n}$ is bounded, which implies that $\beta_n$ cannot converge to zero, hence $\beta_0>0$.  Then, using the continuity of $\Lambda_{\beta,p}$, we can pass to the limit in \rife{eqbetas} and we deduce that $\beta_0$ is  the spectral exponent with $p=p_0$, i.e. $\beta_0= \beta_{S,p_0}$. This proves that $\beta_{S,p}$ is continuous with respect to $p$.
\qeda
\vskip0.5em

We are now ready to prove \rth{ExistII}.
\vskip0.5em

\noindent{\note{\it Proof of \rth{ExistII}}} 
\vskip0.4em
\noindent {\it Step 1: construction of a  solution.} We use similar  ideas as  in the proof of \rth{Exist}, i.e.  a topological degree  argument. On the Banach space
$X =C_0^1(\overline S)$ (endowed with its natural norm) with positive cone $K$, we set 
\bel{B1}\CB(\gw)=-div_g\left((\gb^2\gw^2+|\nabla\gw|^2)^{\frac{p}{2}-1}\nabla \gw\right)+\gb^2(\gb^2\gw^2+|\nabla \gw|^2)^{\frac{p}{2}-1}\gw+ |\gw|^{q-1}\gw
\ee
$$
 \Psi(\gw):={\cal B}^{-1}\left( \beta\left(\la(\beta)+\beta\right) (\gb^2\gw^2+|\nabla\gw|^2)^{\frac{p}{2}-1}\gw_+\right)\,.
$$
Clearly, $\Psi(w)=w$ implies that $w\geq 0$ and solves \rife{absor}. Then, it is enough to prove the existence of a non trivial fixed point for $\Psi$. Observe that, as in \rth{Exist}, $\Psi$ is a continuous compact operator in $X$ thanks to the $C^{1,\alpha}$ estimates for $p$--Laplace operators, and $\Psi(K)\subset K$.

We now wish to compute the degree of $I-\Psi$. First of all we consider, if $R$ is sufficiently large,   ${\rm deg}(I-\Psi, B^+_R, 0)$ where $B^+_R=B_R\cap K$. 
To this purpose, define, for $t\in [0,1]$, $\Psi^*(\gw,t)=t\Psi(\gw)$. Then
 $\Psi^*$ is a compact map on $X\times[0,1]$ and if $\Psi^*(\gw, t) =\gw$, we have 
\be\label{pbt}
 \BA {l}
-div_g\left((\gb^2\gw^2+|\nabla\gw|^2)^{\frac{p}{2}-1}\nabla \gw\right)+ \gb^2(\gb^2\gw^2+|\nabla \gw|^2)^{\frac{p}{2}-1}\gw+  \frac 1{t^{q-(p-1)}} \gw^q \\
\m
\phantom{--------------}
=t^{p-1} \beta\left(\la(\beta)+\beta\right)(\gb^2\gw^2+|\nabla \gw|^2)^{\frac{p}{2}-1}\gw.
 \EA
 \ee
We get, by the maximum principle, 
 $$
 \left\|\frac \gw t\right\|_\infty^{q-(p-1)} \leq  t^{p-1}\beta^{p-1}\left(\la(\beta)+\beta\right)\leq  \beta^{p-1}\left(\la(\beta)+\beta\right)\,.
 $$
Since $t\leq 1$,   we deduce in particular that   $\|\gw\|_\infty$ is bounded independently on $t$. Then, we have  
$$
\frac 1{t^{q-(p-1)}} \gw^q\leq  \left\|\frac \gw t\right\|_\infty^{q-(p-1)}\|\gw\|_\infty^{p-1} \leq C\|\gw\|_\infty^{p-1}\leq C\,.
$$
 \vskip0.4em
Multiplying by $\omega$ we obtain a similar bound  for $\|\gw\|_{W^{1,p}_0(S)}$, and the regularity theory for $p$--Laplace type equations yields a further  estimate on $\| \nabla\gw\|_\infty$. Therefore, we conclude that there exists a constant $M$, independent on $t\in [0,1]$, such that $t\Psi(\gw) =\gw$ implies $\|\gw\|_X\leq M$. 
As a consequence,  if $R$ is sufficiently large we have $t\Psi(\gw) \neq \gw$ on $\partial B_R$. We deduce that ${\rm deg}(I-t\Psi , B^+_R, 0)$ is constant. Therefore     
\be\label{degR}
{\rm deg}(I-\Psi, B^+_R, 0)={\rm deg}(I-t\Psi , B^+_R, 0)={\rm deg}(I, B^+_R, 0)=1\,.
\ee
\smallskip

\noindent Next, we compute ${\rm deg}(I-\Psi, B^+_r, 0)$ for small $r$. We
set
\bel{B2}\
\CB_t(\gw)=-div_g\left((\gb^2\gw^2+|\nabla\gw|^2)^{\frac{p}{2}-1}\nabla \gw\right)+\gb^2(\gb^2\gw^2+|\nabla \gw|^2)^{\frac{p}{2}-1}\gw+ t|\gw|^{q-1}\gw
\ee
and 
$$
 F(\gw,t):={\CB_t}^{-1}\left( \beta\left(\la(\beta)+\beta\right)\gw_+ (\gb^2\gw^2+|\nabla\gw|^2)^{\frac{p}{2}-1} \right)\,.
$$
Again, we have $\Psi(\cdot)=F(\cdot, 1)$. We claim that there exists  a small $r>0$ such that $F(\gw,t)\neq \gw$ for every $t\in [0,1]$ and $\gw\in \partial B_r$. Indeed, reasoning by contradiction, if  this were not true there would exist a nonnegative  sequence $\gw_n$ such that $0\neq \|\gw_n\|\to 0$, and $t_n\in [0,1]$ such that $F(\gw_n,t_n)=\gw_n$, which means that
$$
\BA {l}
-div_g\left((\gb^2\gw_n^2+|\nabla\gw_n|^2)^{\frac{p}{2}-1}\nabla \gw_n\right)+ \gb^2(\gb^2\gw_n^2+|\nabla \gw_n|^2)^{\frac{p}{2}-1}\gw_n+ t_n\gw_n^{q} \\
\m
\phantom{---------------} 
= \beta\left(\la(\beta)+\beta\right)\gw_n (\gb^2\gw_n^2+|\nabla\gw_n|^2)^{\frac{p}{2}-1}
 \EA
$$
Dividing by $\|\gw_n\|^{p-1}$ and letting $n\to \infty$, we find that $\frac{\gw_n}{\|\gw_n\|}$ would converge to some function $\hat w$ such that $\hat w\geq 0$, $\|\hat w\|=1$ and
$$
\BA {l}
-div_g\left((\gb^2\hat \omega^2+|\nabla\hat\omega|^2)^{\frac{p}{2}-1}\nabla \hat \omega\right)+ \gb^2(\gb^2\hat\omega^2+|\nabla \hat\omega|^2)^{\frac{p}{2}-1}\hat \omega \\
\m
\phantom{---------------} 
= \beta\left(\la(\beta)+\beta\right)\hat \omega (\gb^2\hat\omega^2+|\nabla\hat\omega|^2)^{\frac{p}{2}-1}
 \EA
$$
%for some $t_0\in [0,1]$. 
By  \rth{spec} this means that $\la(\beta)=\Lambda_\beta$, which is not possible since $\la(\beta)>\Lambda_\beta$ because $\beta>\beta_S$ (see Remark \ref{moniff}). Therefore, we conclude that
$F(\gw,t)\neq \gw$ for every $t\in [0,1]$ and $\gw\in \partial B_r$ provided $r$ is sufficiently small. We deduce that
${\rm deg}(I-F(\cdot,t), B_r, 0)$ is constant and in particular 
$$
{\rm deg}(I-\Psi, B^+_r, 0)= {\rm deg}(I-F(\cdot,0), B^+_r, 0)\,.
$$
In order to compute this degree, we perform an homotopy acting on $p$ and $\beta$ by setting $p_t=2t+(1-t)p$ and by taking  $\gb_t$ so that 
$t\mapsto \gb_t$ is continuous on $[0,1]$, $\gb_0=\gb$, $\gb_t>\gb_{S,p_t}$ for every $t\in [0,1]$ (where $\gb_{S,p_t}$ is the spectral exponent for $S$ with $p=p_t$)   and $\gb_1>0$ is large enough.
It follows from \rlemma{contp} that  $\gb_{S,p_t}$ is a continuous function of $t$ and remains bounded as $t\in [0,1]$. Therefore,  a similar choice of function $\gb_t$ is possible. In the space $C_0^1(\overline S)$ we define the mapping
%\footnote{to ensure the continuity and compactness properties of $\CC_t$ we need that the $C^{1,\alpha}$ estimates are uniform with respect to $p_t$ and $\beta_t$. Is there any reference ? } 
$\CC_t$ by
\bel{X1}
\CC_t(\gw):=-div_g\left((\gb_t^2\gw^2+|\nabla\gw|^2)^{\frac{p_t}{2}-1}\nabla\gw\right)+\gb_t^2(\gb_t^2\gw^2+|\nabla\gw|^2)^{\frac{p_t}{2}-1}\gw.
\ee
We set
\bel{X2}\tilde F(\gw,t)=\CC_t^{-1}\left(\gb_t(\gl(\gb_t)+\gb_t)(\gb_t^2\gw^2+|\nabla\gw|^2)^{\frac{p_t}{2}-1}\gw\right).
\ee
Combining the Tolksdorf's construction \cite{To1} which shows  the  uniformity with respect to $p_t$ of the $C^{1,\alpha}$ estimates (with $\ga=\ga_t\in (0,1)$), with  the perturbation method of \cite[Th A1]{PoVe}, we obtain that $(\gw,t)\mapsto\tilde F(\gw,t)$ is compact in $C_0^1(\overline S)\ti [0,1]$. Since $\beta_t>\beta_{S,p_t}$,
clearly $I-\tilde F(.,t)$ does not vanish on $\norm\gw_X=r$ for any $r>0$ which implies that 
$$
{\rm deg}(I-\Psi, B^+_r, 0)= {\rm deg}(I-\tilde F(\cdot,0), B^+_r, 0)= {\rm deg}(I-\tilde F(\cdot,1), B^+_r, 0).
$$
But
\bel{X3}I-\tilde F(\cdot,1)=I-\gb_{1}(\gl(\gb_1)+\gb_1)(-\Gd_g+\gb^2_1)^{-1}.
\ee
Since $-\Gd_g$ has only one eigenvalue in $S$ with positive eigenfunction and multiplicity one, choosing $\beta_1$ large in a way that $\la(\beta_1)\beta_1>\la_1(S)$ it follows that 
\medskip

$$
{\rm deg}(I-\tilde F(\cdot,1), B^+_r, 0)=-1={\rm deg}(I-\Psi, B^+_r, 0)\,.
$$
To conclude, since we have
$$
{\rm deg}(I-\Psi, B^+_R\setminus \overline B^+_r, 0)= {\rm deg}(I-\Psi, B^+_R, 0)- {\rm deg}(I-\Psi, B^+_r, 0) \neq 0
$$
we deduce the existence of some $\gw$ such that $r<\|\gw\|<R$ which is a solution of \rife{absor}.

%%%%%%%%%%%%%%%%%%%%%%%%%%%%%%%%%%%%%%%%%%%%%%%%%%%%%%%%%%%%%%%%%STEP2%%%%%%%%%%%%%%%%%%%%%%%%%%%%%%%%%%%%%%%%%%%%%%%%%%%%%%%%%%%%%%%%%%%%%%%%%%%%%%%%%%%%%%%%%%%%%%%%%
\medskip

 \noindent {\it Step 2: uniqueness.}  If $\gw$ is any positive solution, then $\gb^2\gw^2+|\nabla\gw^2|$ is   positive in $\overline S$. This is obvious in $S$ and it is a consequence of Hopf boundary lemma on $\prt S$. 
 Let $\overline \gw$  and $\gw$ be two positive solutions. Either the two functions are ordered or their graphs intersect. Since all the solutions are positive in $S$ and satisfy Hopf boundary lemma,   we can define
 $$\gth:=\inf\{s\geq 1:s\gw\geq \overline\gw\},$$ 
and denote $\gw^*:=\gth\gw$. Either the graphs of $\overline\gw$ and $\gw^*:=\gth\gw$ are tangent at some interior point $\ga\in S$, or $\gw^*>\overline\gw$ in $S$ and there exists $\ga\in \prt S$ such that 
 $\overline\gw_\gn(\ga) =\gw^*_\gn(\ga) <0$. We put $w=\overline\gw-\gw^*$ and use local coordinates $(\gs_1,...,\gs_d)$ on $M$ near $\alpha$. We denote by $g=(g_{ij})$  the
metric tensor on $M$ and $g^{jk}$ its contravariant components. Then, for any $\varphi\in C^1(S)$,
$${\abs {\nabla\varphi }^{2}}=
\sum_{j,k}g^{jk}\frac {\partial \varphi}{\partial \sigma_{j}}\frac {\partial 
\varphi}{\partial \sigma_{k}}=\langle\nabla\varphi,\nabla\varphi\rangle_g.
$$
If $X=(X^1,...X^d)\in C^1(TM)$ is a vector field, if we lower indices by setting 
$\displaystyle {X^\ell=\sum_{i}g^{\ell i}X_{i}}$, then   
$$div_gX=\frac {1}{\sqrt {\abs g}} \sum_{\ell}
\frac {\partial}{\partial \sigma_{\ell}}\left(\sqrt {\abs g}X^\ell\right)
=\frac {1}{\sqrt {\abs g}} \sum_{\ell,i}\frac {\partial}{\partial \sigma_{\ell}}\left(\sqrt {\abs g}g^{\ell i}X_{i}\right).$$
By the mean value theorem applied to
$$t\mapsto \Phi(t)=\left(\beta^{2}{(\gw^*+tw)}^2+\abs 
{\nabla (\gw^*+tw)}^2\right)^{(\frac{p}{2}-1)}(\gw^*+tw)\qquad t=0,1,
$$
we have, for some $t\in (0,1)$,
 \begin {eqnarray*}
 (\beta^{2}\overline\gw^2+\abs 
{\nabla\overline\gw}^2)^{(\frac{p}{2}-1)} \overline\gw-
(\beta^{2}{\gw^*}^2+\abs 
{\nabla\gw^*}^2)^{(\frac{p}{2}-1)}\gw^*
=\sum_{j}a_{j}\frac {\partial w}{\partial 
\sigma_{j}}+bw,
\end {eqnarray*}
where 
$$b=\left(\beta^{2}{(\gw^*+tw)}^2+\abs {\nabla (\gw^*+tw)}^2\right)^{(\frac{p}{2}-2)}
\left((p-1)\gb^2(\gw^*+tw)^2+\abs{\nabla (\gw^*+tw)}^2\right)
$$
and
$$a_j=(p-2)\left(\beta^{2}{(\gw^*+tw)}^2+\abs 
{\nabla (\gw^*+tw)}^2\right)^{(\frac{p}{2}-2)}(\gw^*+tw)\sum_{k}g^{jk} \myfrac{\prt(\gw^*+tw)}{\prt\gs_k}$$
Considering now
$$t\mapsto \Phi_i(t)=\left(\beta^{2}{(\gw^*+tw)}^2+\abs 
{\nabla (\gw^*+tw)}^2\right)^{(\frac{p}{2}-1)}\myfrac{\prt (\gw^*+tw)}{\prt\gs_i}    \qquad t=0,1,
$$
we see that there exists some $t_i\in (0,1)$ such that
 \begin {eqnarray*}
 (\beta^{2}\overline\gw^2+\abs 
{\nabla\overline\gw}^2)^{(\frac{p}{2}-1)} \myfrac{\prt \overline\gw}{\prt\gs_i}-
 (\beta^{2}{\gw^*}^2+\abs 
{\nabla\gw^*}^2)^{(\frac{p}{2}-1)} \myfrac{\prt \gw^*}{\prt\gs_i}
=\sum_{j}a_{ij}\frac {\partial w}{\partial 
\sigma_{j}}+b_iw,
\end {eqnarray*}
where 
\begin {eqnarray*}
b_i=(p-2)\left(\gb^2(\gw^*+t_i w)^2+{\abs 
{\nabla(\gw^*+t_iw)}}^2\right)^{(\frac{p}{2}-2)}\beta^2(\gw^*+t_iw)\frac {\partial 
(\gw^*+t_iw)}{\partial\sigma_{i}}
\end {eqnarray*}
and
\begin {eqnarray*}
a_{ij}=(p-2)\left(\gb^2(\gw^*+t_iw)^2+{\abs 
{\nabla(\gw^*+t_iw)}}^2\right)^{(\frac{p}{2}-2)}\frac {\partial 
(\gw^*+t_iw)}{\partial\sigma_{i}}\sum_{k}g^{jk}\frac {\partial 
(\gw^*+t_iw)}{\partial\sigma_{k}}\\
+\delta_{i}^j\left(\gb^2(\gw^*+t_iw)^2+{\abs 
{\nabla(\gw^*+t_iw)}}^2\right)^{(\frac{p}{2}-1)}.
\end {eqnarray*}
Set $P=\gw^*(\ga)=\overline\gw(\ga)$ and $Q=\nabla\gw^*(\ga)=\nabla\overline\gw(\ga)$. Then $P^2+\abs Q^2>0$ and 
$$b_i(\ga)=
(p-2)\left(\beta^2P^2+\abs {Q}^2\right)^{(\frac{p}{2}-2)}\beta^2\, PQ_{i},
$$
and
$$a_{ij}(\ga)=\left(\beta^2P^2+{\abs Q}^2\right)^{\frac{p}{2}-2}
\left(\delta_{i}^j(\beta^2P^2+{\abs Q}^2)+(p-2) Q_{i}\sum_{k}g^{jk}Q_{k}\right).
$$
Because $\gw^*$ is a supersolution for $(\ref{absor})$, the function $w$ satisfies
\begin {equation}\label{W}
-\frac {1}{\sqrt {\abs g}}\sum_{\ell,j}\frac {\partial}{\partial 
\sigma_{\ell}}\left(A_{j\ell}\frac {\partial w}{\partial \sigma_{j}}\right)
+\sum_{i}C_{i}\frac {\partial w}{\partial 
\sigma_{i}} +Dw\leq 0
\end {equation}
where the $C_{i}$ and $D$ are continuous functions and 
$$A_{j\ell}=\sqrt{\abs g}\sum_{i}g^{\ell i} a_{ij}.
$$
The matrix $(a_{ij})(a)$ is symmetric definite and positive since it is the Hessian of 
$$x=(x_1,...,x_d)=\myfrac{1}{p}(P^2+|x|^2)^{\frac{p}{2}}=\myfrac{1}{p}\left(P^2+\sum_{j,k}g^{jk}x_jx_k\right)^{\frac{p}{2}}
$$
Therefore the matrix $(A_{j\ell})$ keeps the same property in a neighborhood of $a$. Since $w$ is nonpositive and vanishes at some $a\in S$ or $w<0$ and $w_\gn=0$ at some boundary point, it  follows from the strong maximum principle or Hopf boundary lemma (see \cite{PW}) that $w\equiv 0$, i.e. $\theta \omega=\overline\gw$. This implies that actually $\gth=1$ and $\omega=\overline\gw$. 
 \qeda
 %%%%%%%%%%%%%%%%%%%%%%%%%%%%%%%%%%%%%%%%%%%%%%%%%%%%%%%%%%%%%%%%%%%%%%%%%%%%%%%%%%%%%%%%%%%%%%%%%%%%%%%%%%%%%%%%%%%%%%%%%%%%%%%%%%%APPENDIX%%%%%%%%%%%%%%%%%%%%%%%%%%%%%%%%%%%%%%%%%%%%%%%%%%%%%%%%%%%%%%%%%%%%%%%%%%%%%%%%%%%%%%%%%%%%%%%%%%%%%%%%%%%%%%%%%%%%%%%%%%%%%%%%%%%%%%%%%%%%%%%%%%%%%%%%%%%%%%
\section{Appendix}
\setcounter{equation}{0}

We prove here the following result
\bth{isol} Let $S$ be a subdomain of a complete d-dimensional Riemannian manifold $(M,g)$. If $\gb>0$ and $p>1$, the first eigenvalue $\gl_{1,\gb}$ of the operator $\gw\mapsto-div ((\gb^2\gw^2+|\nabla\gw|^2)^{\frac{p}{2}-1}\nabla\gw)+\gb^2\gw(\gb^2\gw^2+|\nabla \gw|^2)^{\frac{p}{2}-1}$ in $W^{1,p}_0(S)$ is isolated. Furthermore any corresponding eigenfunction has constant sign.
\es

\noindent\Proof The proof is an adaptation of the original one due to Anane and Lindqvist when $\gb=0$. 
We recall that 
\begin{equation}\label{A1}
\gl_{1,\gb}=\inf\left\{\myint{S}{}(\gb^2\gw^2+|\nabla\gw|^2)^{\frac{p}{2}}dv_g:\gw\in W^{1,p}_0(S),
\myint{}{}|\gw|^{p}dv_g=1\right\},
\end{equation}
and that there exists $\gw\in W^{1,p}_0(S)\cap C^{1,\ga}(S)$ such that
\begin{equation}\label{A2}
-div ((\gb^2\gw^2+|\nabla\gw|^2)^{\frac{p}{2}-1}\nabla\gw)+\gb^2\gw(\gb^2\gw^2+|\nabla \gw|^2)^{\frac{p}{2}-1}=\gl_{1,\gb}|\gw|^{p-2}\gw\quad\mbox{in } S.
\end{equation}
The function $|\gw|$ is also a minimizer for $\gl_{1,\gb}$, thus it is a positive solution of $(\ref{A2})$. By Harnack inequality \cite{Se}, for any compact subset $K$ of $S$, there exists $C_K$ such that 
$$\myfrac{|\gw|(\gs_1)}{|\gw|(\gs_2)}\leq C_K\qquad\forall \gs_i\in K,\,i=1,2.
$$
Thus any minimizer $\gw$ must keep a constant sign in $S$.
If $\gl_{1,\gb}$ is not isolated, there exists a decreasing sequence $\{\gm_{n}\}$ of real numbers converging to $\gl_{1,\gb}$ and a sequence of functions $\gw_n  \in W^{1,p}_0(S)$, solutions of 
\begin{equation}\label{A3}
-div (\gb^2\gw_n^2+|\nabla\gw_n|^2)^{\frac{p}{2}-1}\nabla\gw_n)+\gb^2\gw_n(\gb^2\gw_n^2+|\nabla \gw_n|^2)^{\frac{p}{2}-1}=\gm_{n}|\gw_n|^{p-2}\gw_n\quad\mbox{in } S
\end{equation}
such that $\norm{\gw_n}_{L^p(S)}=1$. By standard compactness and regularity results, we can assume that 
$\gw_n\to\overline\gw$ weakly in  $W^{1,p}_0(S)$ and strongly in $L^{p}(S)$. Thus
$$\myint{S}{}(\gb^2\overline\gw^2+|\nabla\overline\gw|^2)^{\frac{p}{2}}dv_g\leq \liminf_{n\to\infty}\myint{S}{}(\gb^2\gw_n^2+|\nabla\gw_n|^2)^{\frac{p}{2}}dv_g=\gl_{1,\gb}
$$
which implies that $\overline\gw$ is an eigenfunction associated with $\gl_{1,\gb}$. 

We observe that $\omega_n$ cannot have constant sign. Indeed, if we had that  $\gw_n$ is positive in $\Gw$, we could proceed as in the proof of  \rth{ExistII}-Step 2;  up to rescaling $\omega_n$, we could assume that $w=\gw-\gw_n$ is nonpositive, is not zero,  and the graphs of $\gw$ and $\gw_n$ are tangent. In that case, using \rife{A2} and \rife{A3}, we see that $w$ satisfies a nondegenerate elliptic equation (as in $(\ref{W})$), and we obtain a contradiction either by the strict maximum principle or by Hopf lemma. Thus, any eigenfunction $\omega_n$ must change sign in $\Gw$.
Set $S^+_n=\{\gs\in S:\gw_n(\gs)>0\}$ and $S^-_n=\{\gs\in S:\gw_n(\gs)<0\}$. Clearly, for $0<\gth<1$, 
$$ \BA {l}\myint{S^{\pm}_n}{}(\gb^2\gw_n^2+|\nabla\gw_n|^2)^{\frac{p}{2}}dv_g\geq (1-\gth)\gb^p\myint{S^{\pm}_n}{}|\gw_n|^pdv_g
+\gth\myint{S^{\pm}_n}{}|\nabla \gw_n|^pdv_g.
\EA$$
It follows from $(\ref{A3})$, multiplying by   $\omega_n^+$, that
$$\myint{S^+_n}{}(\gb^2\gw_n^2+|\nabla\gw_n|^2)^{\frac{p}{2}}dv_g=\gm_n\myint{S^+_n}{}|\gw_n|^pdv_g
$$
hence
$$
\gm_n\myint{S^+_n}{}|\gw_n|^pdv_g
\geq (1-\gth)\gb^p\myint{S^+_n}{}|\gw_n|^pdv_g+\gth\myint{S^+_n}{}|\nabla\gw_n|^pdv_g.
$$
Since for some suitable $q>p$ (for example $q=p^*$ if $p<d$, or any $p<q<\infty$ if $p\geq d$)
$$\myint{S^+_n}{}|\nabla\gw_n|^pdv_g\geq c(p,q)\left(\myint{S^+_n}{}|\gw_n|^qdv_g\right)^{\frac{p}{q}}\geq 
c(p,q)|S^+_n|^{\frac{p-q}{q}}\myint{S^+_n}{}|\gw_n|^pdv_g
$$
we obtain
$$\gm_n\geq (1-\gth)\gb^p+\gth c(p,q)|S^+_n|^{\frac{p-q}{q}}.
$$
Similarly we get, multiplying $(\ref{A3})$ by $\omega_n^-$, that
$$\gm_n\geq (1-\gth)\gb^p+\gth c(p,q)|S^-_n|^{\frac{p-q}{q}}.
$$
It follows that the two sets
$$S^{\pm}=\limsup_{n\to\infty}S^{\pm}_n
$$
have positive measure. Since $\overline\gw\geq 0$ on $S^{+}$ and $\overline \gw\leq 0$ on $S^{-}$, we derive a contradiction with the fact that any eigenfunction corresponding to $\gl_{1,\gb}$ has constant sign. \qeda

%%%%%%%%%%%%%%%%%%%%%%%%%%%%%%%%%%%%%%%%%%%%%%%%%%%%%%%%%%%%%%%%%%%%%%%%%%%%%%%%%%%%%%%%%%%%%%%%%%%%%%%%%%%%%%%%%%%%%%%%%%%%%%%%%%%%%%%%%%NEW-SECTION%%%%%%%%%%%%%%%%%%%%%%%%%%%%%%%%%%%%%%%%%%%%%%%%%%%%%%%%%%%%%%%%%%%%%%%%%%%%%%%%%%%%%%%%%%%%%%%%%%%%%%%%%%%%%%%%%%%%%%%%%%%%%%%%%%%%%%%%%%%%%%%%%%%%%%%%%%%%%%%%%%%%%%%%%%%%%%%%%%%%%%%%%%%%%%%%%%%%%%%%%%%%%%%%%%%%%%%%%%%%

\end{document}